\documentclass[12pt]{article}
\usepackage{amsthm,amsmath,amssymb,amscd,verbatim}
\usepackage{graphics,pifont}
\usepackage{amssymb}
\usepackage{graphicx}
\date{}

\begin{document}
\title{On the number of parameters $c$ for which the point $x=0$ is a superstable periodic point of $f_c(x) = 1 - cx^2$}
\author{Bau-Sen Du \cr
Institute of Mathematics \cr
Academia Sinica \cr
Taipei 10617, Taiwan \cr
dubs@math.sinica.edu.tw \cr}
\maketitle

\begin{abstract}
Let $f_c(x) = 1 - cx^2$ be a one-parameter family of real continuous maps with parameter $c \ge 0$.  For every positive integer $n$, let $N_n$ denote the number of parameters $c$ such that the point $x = 0$ is a (superstable) periodic point of $f_c(x)$ whose least period divides $n$ (in particular, $f_c^n(0) = 0$).  In this note, we find a recursive way to depict how {\it some} of these parameters $c$ appear in the interval $[0, 2]$ and show that $\liminf_{n \to \infty} (\log N_n)/n \ge \log 2$ and this result is generalized to a class of one-parameter families of continuous real-valued maps that includes the family $f_c(x) = 1 - cx^2$.  

\bigskip
\noindent{{\bf Keywords}: Bubbles, point bifurcations, (superstable) periodic points, Implicit function theorem}

\medskip
\noindent{{\bf AMS Subject Classification}: 37E05; 37G15; 58F20}
\end{abstract}

The one-parameter family of logistic maps $q_\lambda(x) = \lambda x(1-x)$ (which is topologically conjugate to the family $f_c(x) = 1 - cx^2$) has been used by Verhulst to model population growth and has many applications in modern mathematics, physics, chemistry, biology, economics and sociology {\bf\cite{aus}}.  It is well-known that {\bf\cite{BC, Dev, Ela}} the periodic points of this family are born through either period-doubling bifurcations or saddle-node (tangent) bifurcations and the computation of the exact bifurcation points of this family is a formidable task {\bf\cite{kk}}.  On the other hand, the (least) periods of the first appearance of these periodic points follow the Sharkovsky ordering {\bf\cite{sh}}.  However, it is not clear how the periods of the {\it latter} periodic points appear except through the symbolic MSS sequences of the superstable periodic points as introduced in {\bf\cite{MSS}}.  In this note, we present a more intuitive and quantitative interpretation of this.  

Let $f_c(x) = 1 - cx^2$ be a one-parameter family of continuous maps from the real line into itself with $c$ as the parameter.  By solving the equation $f_c(x) = x$, we obtain that $x_\pm(c) = (-1 \pm \sqrt{\,4c+1})/(2c)$.  When $c > -1/4, f_c(x)$ has two distinct fixed points and these fixed points are born (through tagent bifurcation) from nowhere at $c = -1/4$.  Note that the fixed points $x_-(c)$ of $f_c(x)$ is stable for $-1/4 < c < 3/4$ because for $c$ in this range we have $|f_c'(x_-)| < 1$.  So, when the fixed points of $f_c(x)$ are born, one of them is stable for a while.

We next solve the equation $f_c^2(x) = x$.  This equation is a polynomial equation of degree 4 whose solutions contain fixed points of $f_c(x)$.  So, the quadratic polynomial $1 - cx^2 - x$ must be a factor of $x - f_c^2(x)$.  By solving $x = f_c^2(x) = 1 - c(1-cx^2)^2 = 1 - c(1 - cx^2 - x + x)^2$, we obtain  $0 = (1-cx^2-x)[c^2x^2 - cx - (c-1)]$.  So, $x_\pm^*(c) =(1\pm \sqrt{4c-3})/(2c)$ are periodic points of $f_c(x)$ with least period 2 which must form a period-2 orbit of $f_c(x)$.  This orbit exists for all $c > 3/4$ and is born at $c= 3/4$ from the fixed point $x = x_-(c)$ right after the stable fixed point $x_-(c)$ of $f_c(x)$ loses its stability.  Furthermore, $|(f_c^2)'(x_\pm^*)| < 1$ when $3/4 < c < 5/4$.  That is, the period 2 orbit $\{x_+^*(c), x_-^*(c) \}$ takes on the stability of the fixed points $x_-(c)$ right after it is born. 

We now want to find the periodic points of $f_c(x)$ with least period 3.  As above, we solve the equation $f_c^3(x) = x$.  After some calculations, we obtain that $f_c^3(x) - x = (1-x-cx^2)h(c,x)$, where $h(c,x) =  c^6x^6 - c^5x^5 + (-3c^5+c^4)x^4+(2c^4-c^3)x^3+(3c^4-c^3+c^2)x^2+(-c^3+2c^2-c)x-c^3+2c^2-c+1$.  So, for any fixed $c$, the real solutions (if any) of $h(c, x) = 0$ will be the periodic points of $f_c(x)$ with least period 3 (in particular, when $c_3 \approx 1.7549$ is the unique positive zero of the polynomial $x^3-2x^2+x-1$, the point $x = 0$ is a period-3 point of $f_{c_3}(x)$).  But $h(c,x)$ is a polynomial in $x$ (with $c$ fixed) of degree 6. It is almost impossible to solve it as we did above for fixed points and period-2 points.  Fortunately, with the help of Implicit Function Theorem, at least we can find when these period-3 orbits are born {\bf\cite{Bech, Du1, Li1, Mir, Saha}} and exist for how long {\bf\cite{Gor}}.  By solving the equations $h(c,x) = 0$ and $\frac {\partial}{\partial x} h(c, x)= 0$ simultaneously, we obtain that $c = 7/4$ and $h(7/4,x) = (1/64)^2[343x^3-98x^2-252x+8)]^2$.  Conversely, if $h(7/4,x) = (1/64)^2[343x^3-98x^2-252x+8)]^2$, then $h(7/4,x)$ has three distinct real zeros and for each of such real zero $x$, we have $h(7/4,x) = 0$ and $\frac {\partial}{\partial x} h(7/4, x)= 0$.  Since there are 3 changes in signs of the coefficients of $h(2, x)$.  So, $h(2, x)=0$ has at least one and hence 6 real solutions.  By Implicit Function Theorem, we can continue each of these 6 solutions further from $c = 2$ as long as $\frac {\partial}{\partial x} h(c, x) \ne 0$  and this inequality holds as long as $c > 7/4$.   Again, by Implicit Function Theorem, since $h(0, x) (\equiv 1)$ has no real zeros, $h(c,x) = 0$ has no real solutions for any $0 \le c < 7/4$.  Therefore, the period-3 orbits of $f_c(x)$ are born at $c = 7/4$ and exist for all $c \ge 7/4$.

In theory, we can proceed as above to find the bifurcations of periodic orbits of periods $m \ge 4$.  However, in practice, it becomes more and more difficult as the degree of $f_c^n(x)$ which is $2^n$ grows exponentially fast.  Surprisingly, by extending an idea of Lanford {\bf{\cite{lan}}}, we can show that the number of parameters $c$ such that the point $x = 0$ is a periodic point of $f_c$ with some period grows exponentially fast with its periods.  Indeed, let $\mathcal R$ denote the real line and let $\mathcal X$ be the class of all continuous maps $\phi_c(x) = \phi(c, x) : [0, 2] \times \mathcal R \longrightarrow \mathcal R$ considered as one-parametr families of continuous maps from $\mathcal R$ into itself with parameter $c \in [0, 2]$ such that 

\begin{itemize}
\item[(a)]
$\phi_c(0) > 0$ for all $0 \le c \le 2$; 

\item[(b)]
there exists a {\it smallest} integer $r \ge 2$ such that $\phi_{\hat c}^r(0) = 0$ for some parameter $0 \le \hat c < 2$; and 
 
\item[(c)]
$\phi_2^n(0) < 0$ for all integers $n \ge 2$.    
\end{itemize}
It is clear that the family $f_c(x) = 1 - cx^2$ is a member of $\mathcal X$ since when $r = 2$ and $c = 1$, $\{0, 1\}$ is a period-2 orbit of $f_1(x) = 1 - x^2$ and $f_2^n(-1) = -1 < 0$ for all integers $n \ge 2$. Now let $\Phi_n(c) = \phi_c^n(0)$ for all integers $n \ge 1$ and all $0 \le c \le 2$.  Then, for each integer $n \ge 1$, $\Phi_n(c)$ is a continuous map from $[0, 2]$ into $\mathcal R$ and the solutions of $\Phi_n(c) = 0$ are parameters $c$ for which the point $x = 0$ is a periodic point of $\phi_c(x)$ whose least period divides $n$.  For $\Phi_n(c)$, we have the following 4 properties:
\begin{itemize}

\item[(1)] $\Phi_1(c) = \phi_c(0) > 0$ for all real numbers $0 \le c \le 2$, $r$ is the {\it smallest} integer $n \ge 2$ such that $\Phi_n(c) = 0$ has a solution in $[0, 2)$ and $\Phi_n(2) < 0$ for all integers $n \ge 2$.  

\item[(2)] If $\Phi_{n-1}(\hat c) = 0$ for some integer $n \ge r+1$, then $\Phi_n(\hat c) = \phi_{\hat c}^n(0) = \phi_{\hat c}(\phi_{\hat c}^{n-1}) = \phi_{\hat c}(\Phi_{n-1}(\hat c)) = \phi_{\hat c}(0) > 0$ and so, by (1), $\Phi_n(c) = 0$ has a solution in the interval $(\hat c, 2)$.

\item[(3)] For each integer $n \ge r+1$, the largest solution of $\Phi_n(c) = 0$ is larger than any solution of $\Phi_{n-1}(c) = 0$ and hence than any solution of $\Phi_m(c) = 0$ for any $r \le m < n$.  So, if $c_n^*$ is the largest solution of $\Phi_n(c) = 0$ in the interval $(0, 2)$, then the point $x = 0$ is a periodic point of $\phi_{c_n^*}(x)$ with least period $n$. 

{\it Proof.}
If $c^*$ is the largest solution of $\Phi_{n-1}(c) = 0$, then since, by (2), $\Phi_n(c^*) > 0$ and by (1), $\Phi_n(2) < 0$, we obtain that $\Phi_n(c) = 0$ has a solution between $c^*$ and 2.  The desired result follows accordingly.
\end{itemize}

\begin{itemize}
\item[$(4)_{n,i}$] Let $n \ge 5$ and $i$ be fixed integers such that $2 \le i \le n-2$.  If $\Phi_{n-1}(c_1) = 0$, $\Phi_{n-i}(c_2) = 0$, and $c_2$ {\it is larger than any solution} of $\Phi_i(c) = 0$, then there is a solution of $\Phi_n(c) = 0$ between $c_1$ and $c_2$. 

{\it Proof.}
By (2), $\Phi_n(c_1) > 0$ and by definition, we have $\Phi_n(c_2) = \phi_{c_2}^n(0) = \phi_{c_2}^i(\phi_{c_2}^{n-i}(0)) = \phi_{c_2}^i(\Phi_{n-i}(c_2)) = \phi_{c_2}^i(0) = \Phi_i(c_2)$.  Since $c_2$ is larger than any solution of $\Phi_i(c) = 0$ and $\Phi_i(2) < 0$, we see that $\Phi_n(c_2) = \Phi_i(c_2) < 0$.  Thus, there is a solution of  $\Phi_n(c) = 0$ between $c_1$ and $c_2$. 
\end{itemize}

Now since $r$ is the {\it smallest} integer $n \ge 2$ such that the equation $\Phi_n(c) = 0$ has a solution in the interval $[0, 2]$, if $r \ge 3$ then the equation $\Phi_i(c) = 0$ has no solutions in $[0, 2]$ for each integer $2 \le i < r$ and so, since $\Phi_n(2) < 0$ for all integer $n \ge 2$ by (1), we have $\Phi_i(c) < 0$ for each integer $2 \le i < r$ and all $0 \le c \le 2$.  This fact will be used below.  Let $s = \max\{3, r\}$.  For each integer $k \ge s$, let $c_k^*$ denote the largest solution of $\Phi_k(c) = 0$.  Then, by (3), $0 \le c_r^* < c_{r+1}^* < \cdots < c_k^* < c_{k+1}^* < \cdots < 2$.  For any integer $n \ge k+2$, the equation $\Phi_n(c) = 0$ may have {\it more than} one solution in $[c_k^*, c_{k+1}^*]$.  {\it In the sequel}, let $c_n$ denote {\it any} solution of them.  To distinquish them later, we shall let $[a : b]$ denote the closed interval with $a$ and $b$ as endpoints, where $a$ and $b$ are distinct real numbers.  When there is no confusion occurs, we shall use the same $c_n$ to denote several distinct solutions of $\Phi_n(c) = 0$ in the real line.  we now describe how to apply $(4)_{n,i}$ successively, in the interval $[c_k^*, c_{k+1}^*]$, to obtain solutions of $\Phi_n(c) = 0$ with $n \ge k+2$.  The procedures are as follows:
\begin{itemize}
\item[(i)]
We start with the interval $[c_k, c_{k+1}] = [c_k^*, c_{k+1}^*]$.

\item[(ii)]
The $1^{st}$ step is to apply $(4)_{k+2,2}$ to $[c_k, c_{k+1}]$ in (i) to obtain {\it one} $c_{k+2}$ between $c_k$ and $c_{k+1}$ and so obtain the 2 intervals $[c_k, c_{k+2}]$ and $[c_{k+2}, c_{k+1}]$.  Then the $2^{nd}$ step is to apply $(4)_{k+3,3}$ to $[c_k, c_{k+2}]$ to obtain the first $c_{k+3}$ in $(c_k, c_{k+2})$ and apply $(4)_{k+3,2}$ to $[c_{k+2}, c_{k+1}]$ to obtain the second $c_{k+3}$ in $(c_{k+2}, c_{k+1})$ and so, we obtain two parameters $c_{k+3}$'s in $(c_k, c_{k+1})$ which, together with the previously obtained one parameter $c_{k+2}$, divide the interval $[c_k, c_{k+1}]$ into 4 subintervals such that $$c_k \quad < \quad \text{first} \,\,\, c_{k+3} \quad < \quad c_{k+2} \quad < \quad \text{second} \,\,\, c_{k+3} \quad < \quad c_{k+1}.$$  Similarly, the $3^{rd}$ step is to apply $(4)_{k+4,i}$ for appropriate $2 \le i \le 4$ to each of these 4 subintervals to obtain 4 parameters $c_{k+4}$'s which, together with the previously obtained $c_{k+3}$'s and $c_{k+2}$, divide the interval $[c_k, c_{k+1}]$ into 8 subintervals such that $$c_k < \text{first} \,\, c_{k+4} < 1^{st} \, c_{k+3} < \text{second} \,\, c_{k+4} < c_{k+2} < \text{third} \,\, c_{k+4} < 2^{nd} \, c_{k+3} < \text{fourth} \,\, c_{k+4} < c_{k+1}.$$We proceed in this manner indefinitely to obtain, at the $i^{th}$ step, $i \ge 1$, several parameters $c_{k+1+i}$'s which are interspersed with parameters $c_j$'s with {\it smaller} subscripts $k \le j \le k+i$ and note that each $c_{k+1+i}$ is {\it ajacent} to a $c_{k+i}$ on the one side and to a $c_j$ with $k \le j < k+i$ on the other.

\item[(iii)]
For any two ajacent parameters $c_i$ and $c_j$ in $[c_k^*, c_{k+1}^*]$ with $k \le i < j \le 2k-2$ (and so, $2 \le j+1-i \le k-1$), we apply $(4)_{j+1, j+1-i}$ to the interval $[c_i : c_j]$ to obtain one $c_{j+1}$ between $c_i$ and $c_j$.  Consequently, inductively for each $2 \le \ell \le k-1$, we can find one parameter $c_{k+\ell}$ from each of the $2^{\ell -2}$ pairwise disjoint open components formed by the previously obtained $c_{k+s}$'s, $2 \le s \le \ell -1$ in the interval $[c_k^*, c_{k+1}^*]$.  In particular, when $\ell = k-1$, we obtain $2^{k-3}$ parameters $c_{2k-1}$'s and $2^{k-2}$ such intervals $[c_{i} : c_{2k-1}]$, $k \le i \le 2k-2$, with mutually disjoint interiors.  Now we apply appropriate $(4)_{2k,2k-i}, k+1 \le i \le 2k-2$, to each of the $2^{k-2}$ previously obtained intervals to obtain $2^{k-2}-1$ parameters $c_{2k}$'s (here the minus 1 is added because we do not have a parameter $c_{2k}$ which is right {\it next} to $c_k$).  In summary, so far, we have obtained, for each integer $k \le i \le 2k$, $a_{k,i}$ parameters $c_{k+i}$ such that they are interspersed in a pattern similar to those depicted in (ii) above, where $a_{k,k} = 1 = a_{k,k+1}, a_{k,k+j} = 2^{j-2}, 2 \le j \le k-1$, and $a_{k,2k} = 2^{k-2} - 1$.  

\item[(iv)]
For each integer $m > 2k$, let $a_{k,m} = \sum_{1 \le i \le k-1} a_{k, m-i}$.  Then it follows from {\bf\cite{w}} that $\lim_{m \to \infty} \log a_{k,m} = \log \alpha_k$, where $\alpha_k$ is the (unique) positive (and largest in absolute value) zero of the polynomial $x^{k-1} - \sum_{i=0}^{k-2} x^i$ and $\lim_{k \to \infty} \alpha_k = 2$.  We now continue the above procedures by applying $(4)_{n,i}$ with appropriate $n$ and $i$ to each appropriate interval $[c_s, c_t]$ as we do in (ii) and (iii) above (see Figure 1).  We want to show that, for each integer $m > 2k$,  there are at least $a_{k,m}$ distinct parameters $c_m$'s in $[c_k^*, c_{k+1}^*]$ which are obtained in this way.  Indeed, given $m > 2k$, let $c_{m-k} \,\, (= c_{k+1+(m-2k-1)})$ be {\it any} parameter which is obtained at the $(m-2k-1)^{st}$ step and let $c_{j}$ be a parameter among all parameters $c_s$ with smaller subscripts $k \le s < m-k$ which is {\it next} to $c_{m-k}$ (on either side of $c_{m-k}$).  Recall that at least one such $c_j$ is $c_{m-k-1}$.  If $m-k-j \ge k$, then we ignore this $c_{m-k}$ because no appropriate $(4)_{n,i}$ can be applied to $[c_{m-k} : c_j] \,\, (\subset [c_k^*, c_{k+1}^*])$ to obtain {\it one} $c_{m-k+1}$.  Otherwise, we apply $(4)_{m-k+1,m-k+1-j}$ to $[c_{m-k} : c_j]$ to obtain a parameter $c_{m-k+1}$ between $c_{m-k}$ and $c_j$ and apply $(4)_{n,i}$ with appropriate $i$ to $[c_{m-k} : c_{m-k+1}]$ successively as $n$ increases from $m-k+2$ to $m$ until we obtain several $c_m$'s interspersing with parameters $c_t$ with smaller subscripts $m-k \le t < m$ which are obtained earlier in $[c_{m-k} : c_j]$.  It is clear that there is a parameter $c_m$ which is {\it next} to $c_{m-k}$, but none is {\it next} to $c_j$.  Therefore, in the interval $[c_{m-k} : c_j]$, there are as many parameters $c_m$'s as the sum of all parameters $c_r$'s with integers $r$ in $[m-k+2, m-1]$.  We conclude that, altogether, there are at least $a_{k,m}$ distinct parameters $c_m$'s in $[c_k^*, c_{k+1}^*]$ such that the point $x = 0$ is a periodic point of $\phi_{c_m}(x)$ whose least period divides $m$.  
\end{itemize}

\begin{figure}[ht] 
\begin{center}
\includegraphics[width=17cm,height=13.5cm]{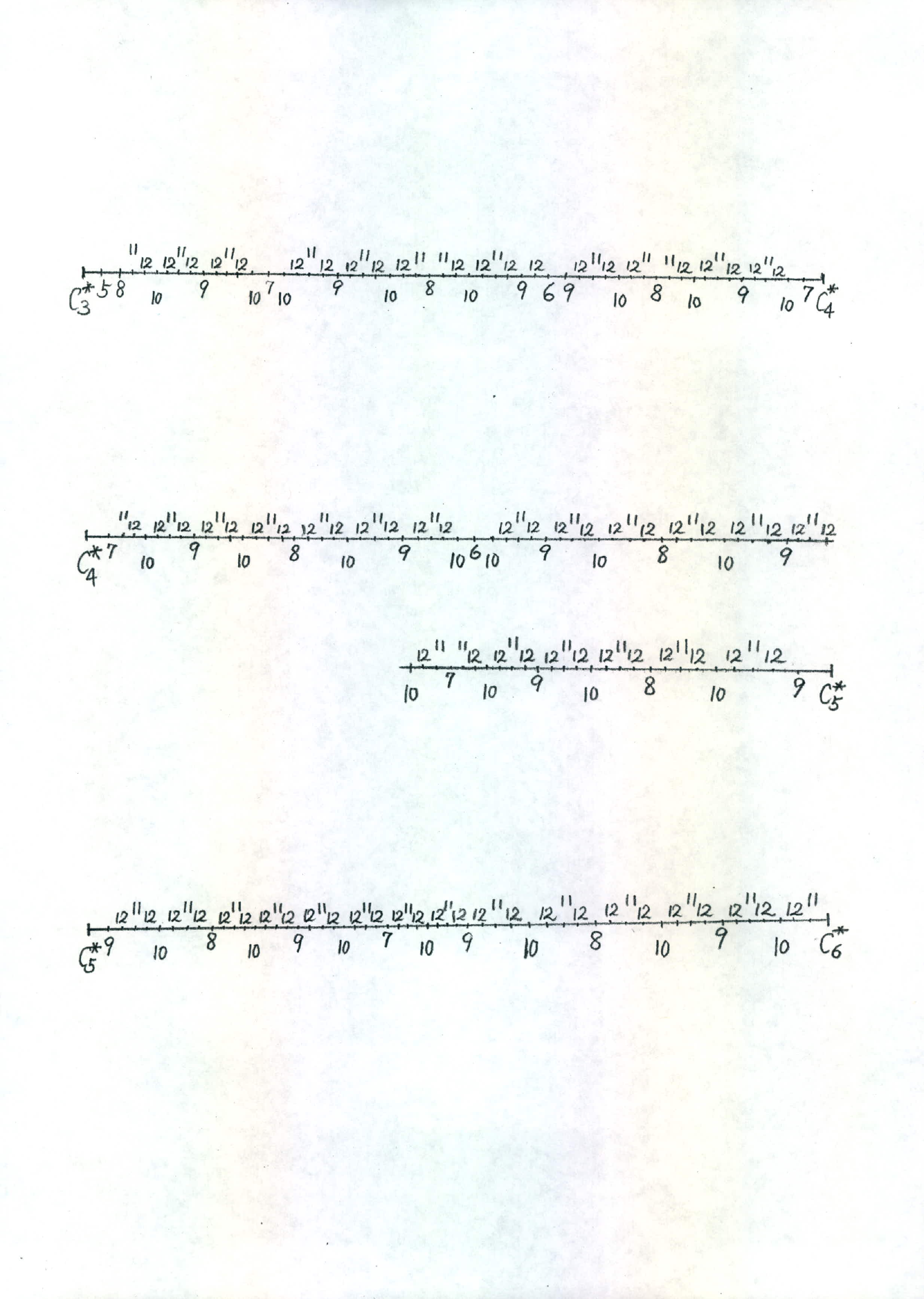} 
\caption{Distribution of the parameters $c_n$'s in the intervals $[c_3^*, c_4^*], [c_4^*, c_5^*]$ and $[c_5^*, c_6^*]$ for $n$ up to 12.  For simplicity, we surpress the letter $c$ and write only the subscript $n$.  So, every positive integer $n$ here is only a symbol representing a parameter $c_n$.}
\end{center}
\end{figure}

In particular, for each integer $n \ge \hat r = \max\{3, r\}$, we have found at least $a_{k,n}$ distinct solutions of $\Phi_n(c) = 0$ in $[c_k^*, c_{k+1}^*]$ for each $\hat r \le k \le n-1$.  So, totally, we have found at least $N_n = \sum_{k=\hat r}^{n-1} a_{k,n}$ distinct solutions of $\Phi_n(c) = 0$ in $[c_{\hat r}^*, c_n^*] \,\, (\subset [0, 2]$).  That is, there are at least $N_n$ distinct parameters $c$ such that $x = 0$ is a periodic point of $\phi_c(x)$ whose least period divides $n$.  Therefore, since $\liminf_{n \to \infty} (\log N_n)/n \ge \sup_{k \ge \hat r} \bigl(\lim_{n \to \infty} (\log a_{k,n})/n\bigr) = \log 2$, we have proved the following result.

\noindent
{\bf Theorem 1.}
{\it Let $\mathcal R$ denote the set of all real numbers and let $\phi(c, x) : [0, 2] \times \mathcal R \longrightarrow \mathcal R$ be a continuous map.  Write $\phi_c(x) = \phi(c, x)$ and consider $\phi_c(x)$ as a one-parameter family of continuous real-valued maps with $c$ as the parameter.  Assume that 
\begin{itemize}
\item[{\rm (a)}]
$\phi_c(0) > 0$ for all $0 \le c \le 2$; 

\item[{\rm (b)}]
there exists a {\it smallest} integer $r \ge 2$ such that $\phi_{\hat c}^r(0) = 0$ for some parameter $0 \le \hat c < 2$; and 

\item[{\rm (c)}]
$\phi_2^n(0) < 0$ for all integers $n \ge 2$.  
\end{itemize}
Let $\hat r = \max\{3, r\}$ and, for each integer $n \ge 2$, let $\Phi_n(c) = \phi_c^n(0)$.  For each integer $\ell \ge r$, let $c_\ell^*$ be the largest solution of the equation $\Phi_\ell(c) = 0$ and, for each integer $k \ge {\hat r}$, let $a_{k,k} = 1 = a_{k, k+1}$, $a_{k, k+i} = 2^{i-2}$, $2 \le i \le k-1$, $a_{k, 2k} = 2^{k-2}-1$, and, for each integer $j > 2k$, let $a_{k, j} = \sum_{1 \le i \le k-1} a_{k, j-i}$.  Also, let $N_{\hat r} = N_{\hat r+1} = 1$ and $N_m = \sum_{k=\hat r}^{m-1} a_{k,m}$ for all integers $m \ge \hat r+2$.  Then the following hold:
\begin{itemize}
\item[{\rm (1)}]
$0 \le c_r^* < c_{r+1}^* < c_{r+2}^* < \cdots < 2$;

\item[{\rm (2)}]
Let $k \ge \hat r$ be a fixed integer.  Then, for each integer $i \ge 0$, there are at least $a_{k,k+i}$ distinct parameters $c_{k+i}$'s in the interval $[c_k^*, c_{k+1}^*]$ such that $\phi_{c_{k+i}}^{k+i}(0) = 0$.  

\item[{\rm (3)}]
For each integer $n \ge \hat r$, there are at least $N_n$ distinct parameters $c$'s such that the point $x = 0$ is a periodic point of $\phi_c(x)$ whose least period divides $n$ and $\liminf_{n \to \infty} (\log N_n)/n \ge \log 2$.
\end{itemize}}

\noindent
{\bf Remarks.}
(1) For the family $f_c(x) = 1 - cx^2$ with parameter $0 \le c \le 2$, it is well-known that $c_2^* = 1$ and $c_3^* \approx 1.7549$ is the unique positive zero of the polynomial $c^3-2c^2+c-1$.  However, it is also well-known that there is an additional parameter $c_4 \approx 1.3107$ which is a zero of the polynomial $-c^7+4c^6-6c^5+6c^4-5c^3+2c^2-c+1 = 1-c[1-c(1-c)^2]^2$.  Since $c_2^* < c_4 < c_3^*$, we can apply $(4)_{n,2}$ to the interval $[c_4, c_3^*]$ successively to obtain additional parameters $c_n, n \ge 5$ which are not counted in Theorem 1 such that $$1 = c_2^* < c_4 < \cdots < c_{2i} < c_{2i+2} < \cdots < c_{2j+1} < c_{2j-1} < \cdots < c_5 < c_3^*.$$
  
(2) For the family $f_c(x) = 1 - cx^2$ with parameter $c \ge 0$, it is generally believed that once a periodic orbit is born it lives forever.  However, this phenomenon does not shared by some other families of polynomials.  For example: let $g_c(x) = \sqrt{7} - c(1 + x^2)$ and $h_c(x) = x^3- 2x + c$.  These two families have phenomena called bubbles (periodic orbits live for a finite time) and point bifurcations (periodic orbits die on birth).  See {\bf\cite{Du3, Li2}} for details.

\end{document}